\documentclass[a4paper,12pt]{article}
\usepackage{epstopdf}
\usepackage{amsmath} 
\usepackage{amsthm}
\usepackage{amssymb}
\usepackage[utf8]{inputenc}
\usepackage[top=2cm, bottom=2cm, left=2cm, right=2cm]{geometry}
\usepackage{enumitem}
\usepackage{tikz}
\usepackage{subfigure}
\usetikzlibrary{patterns}

\theoremstyle{plain}
\newtheorem{tet}{Theorem}
\newtheorem{all}[tet]{Proposition}
\theoremstyle{definition}

\renewcommand{\Re}{\operatorname{Re}}
\renewcommand{\Im}{\operatorname{Im}}

\title{Global stability for the 2-dimensional logistic map}
\author{J\'anos Dud\'as\\
University of Szeged, Hungary}

\begin{document}

\maketitle


\noindent\textbf{Abstract.} For the delayed logistic equation $x_{n+1} = a x_n (a-x_{n-1})$ it is well known that the nontrivial fixed point is locally stable for $1<a\leq 2$, and unstable for $a>2$. We prove that for $1<a\leq 2$ the fixed point is globally stable, in the sense that it is locally stable and attracts all points of $S$, where $S$ contains those $(x_0,x_1)\in \mathbb{R}_+^2$, for which the sequence $\left\lbrace x_n\right\rbrace \subset \mathbb{R}_+$. The proof is a combination of analytical and reliable numerical methods.

\bigskip


\bigskip


\section{Introduction}

One of the most studied nonlinear maps is the logistic map
$$ [0, 1] \ni x \mapsto a x (1-x) \in \mathbb{R},$$
with parameter $a>0$. For $0 < a \leq 1$, it is well known (see e.g. \cite{dev04}) that $x=0$ is the unique fixed point in $[0,1]$, and it is globally stable (i.e. stable and attracts all points in $[0,1]$). For $1 < a \leq 3$, there is a nontrivial fixed point $x_* = 1-\frac{1}{a}$ which is stable and attracts all points in $(0,1)$. At $a=3$ a period doubling (flip) bifurcation takes place, and the fixed point $x_*$ becomes unstable for $a>3$. As $a$ increases, there is a sequence of bifurcation points, and for some larger value of $a$, chaotic behaviour can be shown.

In 1968, Maynard Smith \cite{ms68} considered the "delayed" version
$$
x_{k+1}=a x_k (1-x_{k-1}),
$$
of the logistic difference equation. This is natural in the context of population models: the size of the subsequent generation of the population depends not only on the size in the previous year, but also on the size of the two-year-earlier population.

Introducing $y_k=x_{k+1}$, the second order difference equation is equivalent to
\begin{equation*}
(x_{n+1},y_{n+1}) = F_a(x_n,y_n)
\end{equation*}
with
\begin{equation}  \label{falap}
F_a(x,y) = (y, a y (1-x)).
\end{equation}
We study the map $F_a$ for those $(x,y) \in \mathbb{R}^2_+ = [0,\infty)\times[0,\infty)$ for which all iterates of $F_a$ remain in $\mathbb{R}^2_+$, i.e. $F_a^k(x,y)\in \mathbb{R}^2$, for every $k\in \mathbb{N}$. Here $F_a^k$ denotes the $k$-fold iteration of $F_a$, i.e. $F_a^0 = id$, $F_a^k = F_a\left( F_a^{k-1} \right),\ k\in \mathbb{N}$. As we will see, for $0<a\leq 2$ the set 
\begin{equation*}
S_0=\{(x,y)\in \mathbb{R}^2:\ 0\leq x\leq 1; \quad 0\leq y\leq 1; \quad a y (1-x)\leq 1\}
\end{equation*}
is invariant under $F_a$, that is $F_a(S_0) \subseteq S_0$.

For $0 < a \leq 1$, we have $S_0 = [0,1] \times [0,1]$, the only fixed point in $S_0$ is $(0,0)$, which is locally stable and $F_a^k(x,y) \to (0,0)$ as $k \to \infty$. For $a > 1$, the nontrivial fixed point $(A,A) \in S_0$ with $A = 1-\frac{1}{a}$ appears, which is locally asymptotically stable for $a \in (1,2)$, and it is unstable for $a > 2$. A Neimark--Sacker bifurcation takes place at $a = 2$ (see e.g. in \cite{kuzn} Example 4.3) and for $a > 2$, $a$ is close to $2$ a stable invariant curve appears. As we increase $a$, the size of the invariant curve is getting larger; at about $a=2.27$, the curve touches the $x$-axis, and complicated dynamics occurs. For profound numerical studies, see \cite{aronson1982}, \cite{nsbif}.

The aim of this paper is to show that, for $1 < a \leq 2$, the nontrivial fixed point $(A,A)$ is globally stable in the sense that $(A,A)$ is locally stable, and for each $(x,y)$ in
\begin{equation*}
S = S_a = \{(x,y)\in \mathbb{R}^2: \ 0 \leq x < 1; \quad 0 < y < 1; \quad a y (1-x) < 1\},
\end{equation*}
$F_a^k(x,y) \to (A,A)$ as $k \to \infty$. Consequently, the local stability of $(A,A)$ implies its global stability. For similar results on the global stability of other delayed difference equations, the reader is referred for example to \cite{lizgas}, \cite{liztt}, \cite{ladas3} or \cite{bgk}.

We emphasize that we prove the stability even in the critical parameter value $a=2$. However, we do not consider the case $a>2$. According to numerical studies (\cite{aronson1982}, \cite{nsbif}), the invariant curve is globally stable for parameter $a>2$ close to the critical value $2$.

For $a\leq 2$ the proof of the global stability is a combination of analytical and computer-aided tools. It is based on the method in \cite{bgk} and \cite{bg}. We elaborate the analytical part such that it can be easily applied to similar models. Furthermore, a quite important aim is to have a clear picture of the method in order to be able to prove similar results for higher dimensional models, for example the 3-dimensional logistic map $x_{n+1} = a x_n (1- x_{n-2})$, where further difficulties arise.

With analytical tools we construct an attracting neighbourhood $\mathcal{N}$ around the nontrivial fixed point $(A,A)$. Then we show that every $(x_0, y_0)\in S\setminus \mathcal{N}$ will eventually step into $\mathcal{N}$, that is, there exist an $n\in \mathbb{N}$, such that $(x_n, y_n)\in \mathcal{N}$, where $(x_n, y_n) = F_a^n(x_0, y_0)$. So these points are also in the region of attraction of the fixed point $(A,A)$. We use computer, applying reliable numerical methods, to show the second step. In this context, reliable means, we use interval arithmetic tools to control every occurring numerical error, consequently, the method is suitable to prove mathematical statements. (See e.g. \cite{tuc})

In section \ref{seclin}, for smaller parameter values $a$, i.e. for $a\in (1.5,1.95]$ we use the linearised map to construct the attracting neighbourhood $\mathcal{N}$. However, as we will see it later, this neighbourhood shrinks to the fixed point as $a$ tends to the critical value $2$. Therefore, for parameter values $a$ close to  $2$ this neighbourhood is not big enough for computer use in the second part of the method. Thus we need another approach to construct an attracting neighbourhood $\mathcal{N}$. In section \ref{secnorm}, for these parameter values $a$ close to $2$, we use the normal form of the Neimark--Sacker bifurcation. More precisely, with smooth and invertible maps, we transform map (\ref{falap}) into its normal form, hereby we obtain an attracting neighbourhood $\mathcal{N}$ around the fixed point $(A,A)$, whose size is independent of the parameter $a \in [1.95,2]$.

Since we need the size of the constructed neighbourhood $\mathcal{N}$ for computer use, it is not enough to determine only the lower order terms during the normal form transformation, like we would do in a regular bifurcation analysis. These lower order terms only assure the existence of such a sufficiently small neighbourhood, whose size is not explicitly determined by them. Therefore, it is essential during the transformation to trace the higher order terms and to estimate them as well as possible, in order to obtain a sufficiently big neighbourhood $\mathcal{N}$, since the computer method is more and more compute-intensive and time-consuming, as we get closer to the fixed point.

In section \ref{secpc}, we consider those points, which lie outside the attracting neighbourhood, i.e. the points of $S \setminus \mathcal{N}$. We cover $S$ with finitely many small squares. Considering these squares as vertices of a graph, we introduce a directed graph, which, to a certain extent, describes the behaviour of map (\ref{falap}) on these squares. Therefore we convert the problem of examining infinitely many points into a finite graph problem, which can be handled by computer. To construct the edges of the graph we use reliable numerical methods in order to handle the rounding errors of the computer. We show with the help of this graph that every point from $S$ enters the neighbourhood $\mathcal{N}$ constructed before. With this we will prove our main result:
\begin{tet}\label{fotet}
	For all $a\in (1,2]$ the fixed point $(A,A)$ is locally asymptotically stable, and $\lim\limits_{n\to \infty}F_a^n(x,y)=(A,A)$ for every $(x,y)\in S$, where $A=1-\frac 1 a$.
\end{tet}

\section{Preliminaries}

In this section we study the dynamics of the map (\ref{falap}) for $a > 0$ in the positive quadrant.
Introduce the following sets:
	\begin{align*}
		S & = \left\lbrace (x,y): \ 0 \leq x < 1, \quad 0 < y < 1, \quad a y (1-x) < 1\right\rbrace,\\
		T_0 & = \left\lbrace (x,0): \ x \geq 0 \right\rbrace \cup 
		\left\lbrace (1,y): \ y>0 \right\rbrace \cup \left\lbrace (x,1): \ 0\leq x < 1\right\rbrace \\
		&\quad \cup \left\lbrace (x,y): \ a y (1-x) = 1, \quad 0\leq x \leq 1-\frac{1}{a}\right\rbrace,\\
		T_1 & = \left\lbrace (x,y): \ x\geq 0, \quad 0<y<1, \quad a y(1-x) > 1 \right\rbrace, \\
		T_2 & = \left\lbrace (x,y): \ 0\leq x < 1, \quad y>1 \right\rbrace, \\
		T_3 & = \left\lbrace (x,y): \ x>1, \quad y>0 \right\rbrace.	
	\end{align*}

Clearly, $\mathbb{R}_+^2 = S \cup T_0 \cup T_1 \cup T_2 \cup T_3$, furthermore, for $0<a\leq 1$, $S=[0,1)\times(0,1)$ and $T_1 = \emptyset$.

\begin{center}
	\begin{tikzpicture}[scale=3]
	\draw [->] (-.1,0) -- (1.7,0) node[right] {$x$};
	\draw [->] (0,-.1) -- (0,1.7) node[above] {$y$};
	\draw[] (1,-.01) -- (1,.01);
	\path (1,-.01) node[below] {$1$};
	\draw (-.01,1) -- (.01,1);
	\path (-.01,1) node[left] {$1$};
	\draw[very thick] (0,0) -- (1.69,0);
	\draw[very thick] (1,0) -- (1,1.7);
	\draw[very thick] (0,1) -- (1,1);
	\draw[very thick, domain=0:.5] plot (\x, {1/(2-2*\x)});
	\path (.5,.5) node {$S$};
	\path (.2,.82) node {$T_1$};
	\path (.5,1.2) node {$T_2$};
	\path (1.3,.7) node {$T_3$};
	\path (1.3,1.2) node {$T_0$};
	\draw [->] (1.21,1.16) -- (1.01,1.05);
	\end{tikzpicture}
\end{center}

\begin{all}
	For all $a>0$, we have
	$$
	F_a^4(T_0) = \left\lbrace (0,0) \right\rbrace, \quad F_a(T_1) \subset T_2, \quad F_a(T_2) \subset T_3, \quad F_a(T_3)\cap \mathbb{R}^2_+ = \emptyset.
	$$
	and furthermore, if $a\in (0,2]$ then
	$
	F_a(S) \subset S.
	$
\end{all}

\noindent
\textit{Proof. }
From the definition of $F_a$ it is obvious that $F_a^4(T_0) = \left\lbrace (0,0) \right\rbrace$. It is also straightforward to check the relations $F_a(T_1) \subset T_2, \ F_a(T_2) \subset T_3$ and $F_a(T_3)\cap \mathbb{R}^2_+ = \emptyset$.

If $0 < a \leq 2$ and $(x,y)\in S$ then $0<y<1$, $0<a y (1-x) < 1$ and
\begin{equation*}
	a^2 y (1-x)(1-y) \leq 4 (1-x)\max_{0\leq y\leq 1}y(1-y)\leq 1.
\end{equation*}
Therefore $F_a(S)\subset S$.
\qed

\bigskip{}

Consequently, in the rest of the paper we can assume $(x,y)\in S$. For small $a$, the dynamics in $S$ is quite simple. The following statement easily follows from the monotonicity of $\left\lbrace x_n \right\rbrace_{n=0}^\infty$ for $0 < a \leq 1$.
\begin{all}
	If $0 < a \leq 1$, then for all $(x,y) \in [0,1]\times[0,1]$,
	$$
	F^k_a(x,y) \to (0,0) \text{ as \ } k\to \infty.
	$$
\end{all}

For $1<a\leq 2$, we divide the positive quadrant into four subsets with lines $x=A$, $y=A$, and introduce the following sets:
\begin{align*}
	S_1 & =\{(x,y)\in S :\ x\leq A, \quad y<A \}, \\
	S_2 & =\{(x,y)\in S :\ x< A, \quad A\leq y \}, \\
	S_3 & =\{(x,y)\in S :\ A\leq x, \quad A<y\}, \\
	S_4 & =\{(x,y)\in S :\ A<x, \quad y\leq A\}.
\end{align*}

Clearly $S = \bigcup_{i = 1}^4 S_i \cup \left\lbrace (A,A) \right\rbrace $. Introduce the notation $t_n = (x_n,y_n) = F_a^n(x_0,y_0)$. 
\begin{all} \label{kor}
	For every $t_0 = (x_0, y_0)\in S$ and $a\in (1,2]$ the sequence $\{t_n\}_{n=0}^{\infty}$, defined by $t_{n+1} = (x_{n+1},y_{n+1}) = F_a(x_n,y_n)$ fulfils one of the following cases:
	\begin{enumerate}[label=(\alph*)]
		\item $\lim\limits_{n\to \infty} t_n=(A,A)$,
		\item the sequence $\{t_n\}_{n=0}^{\infty}$ goes around the fixed point along the cycle $S_1 \to S_2 \to S_3 \to S_4 \to S_1$; and in the course of one cycle, there can be more than one elements of the sequence in both $S_1$ and $S_3$, but the number of these elements are finite.	
	\end{enumerate}
\end{all}


\begin{figure}
	\begin{center}
		\begin{tikzpicture}[scale=6]
		\draw  (-.1,0) -- (0,0);
		\draw [thick, dashed] (0,0) -- (1,0);
		\draw [->] (1,0) -- (1.1,0) node[right] {$x$};
		\draw (0,-.1) -- (0,0);
		\draw [thick] (0,0) -- (0,4/7);
		\draw [->] (0,3/7) -- (0,1.1) node[above] {$y$};
		\draw[] (3/7,-.01) -- (3/7,.01);
		\path (3/7,-.01) node[below] {$A$};
		\draw[] (1,-.01) -- (1,.01);
		\path (1,-.01) node[below] {$1$};
		\draw[] (-.01,3/7) -- (.01,3/7);
		\path (-.01,3/7) node[left] {$A$};
		\draw (-.01,1) -- (.01,1);
		\path (-.01,1) node[left] {$1$};
		\draw[thick, dashed] (1,0) -- (1,1);
		\draw[thick, dashed] (3/7,1) -- (1,1);
		\draw (0,1) -- (3/7,1);
		\draw[thick, dashed] (0,.426) -- (3/7,.426);
		\draw[thick] (0,.43) -- (3/7,.43);
		\draw[thick, dashed] (3/7,.43) -- (1,.43);
		\draw[thick] (3/7,.426) -- (1,.426);
		\draw[thick, dashed] (.426,3/7) -- (.426,1);
		\draw[thick] (.43,3/7) -- (.43,1);
		\draw[thick, dashed] (.43,0) -- (.43,3/7);
		\draw[thick] (.426,0) -- (.426,3/7);
		\draw[thick,dashed, domain=0:3/7] plot (\x, {4/(7-7*\x)});
		\draw[gray, thick,dotted, domain=3/7:.5] plot (\x, {4/(7-7*\x)}) node[right] {$ay(1-x)=1$};
		\path (.3,.25) node {$S_1$};
		\path (.3,.55) node {$S_2$};
		\path (.75,.75) node {$S_3$};
		\path (.75,.25) node {$S_4$};
		\draw [gray, ->] (.25,.27) arc [radius=.07, start angle=90, end angle= 360];
		\draw [gray, ->] (.73,.8) arc [radius=.07, start angle=180, end angle= -90];
		\draw [gray, ->] (.3,.3) -- (.3,.5);
		\draw [gray, ->] (.35,.57) -- (.7,.75);
		\draw [gray, ->] (.75,.7) -- (.75,.3);
		\draw [gray, ->] (.7,.25) -- (.35,.25);
		\end{tikzpicture}
	\end{center}
	\caption{The dynamics in $S$}
\end{figure}
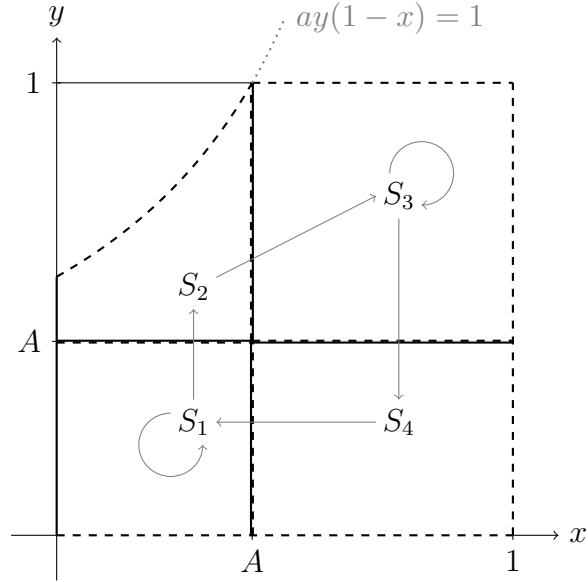

\noindent
\textit{Proof. } The transitions between the aforementioned subsets are the following:

\begin{enumerate}
	\item [--] For $t_0 \in S_1$ we obtain $x_1<A$, therefore $\pmb{S_1 \to \{S_1, S_2\}}$. That is $t_1 \in S_1$ or $t_1 \in S_2$.
	\item [--] For $t_0 \in S_2$: \ $x_1\geq A$ and $y_1=a y_0 (1-x_0)>aA(1-A)=A$, so $\pmb{S_2\to S_3}$.
	\item [--] For $t_0 \in S_3$:  \ $x_1>A$, so $\pmb{S_3 \to \{S_3,S_4\}}$.
	\item [--] For $t_0 \in S_4$: \ $x_1\leq A$ and $y_1=a y_0 (1-x_0)>aA(1-A)=A$, so $\pmb{S_4 \to S_1}$.
\end{enumerate}

We obtain there is a cycle $S_1 \to S_2 \to S_3 \to S_4 \to S_1$. But during a cycle the points of the sequence can spend more time in $S_1$ or $S_3$, possibly, the sequence can stay in $S_1$ or $S_3$ forever. We only need to show that, if a sequence gets stuck in $S_1$ or $S_3$, then it converges to the fixed point $(A,A)$.

Notice that $y_{n+1} \geq y_n$ as long as $x_n\leq A$, and similarly, $x_n\geq A$ implies $y_{n+1} \leq y_n$. According to this, as long as $t_n\in S_1\cup S_2$, the sequence $\{y_n\}=\{x_{n+1}\}$ increases, until the sequence $\{t_n\}$ steps into $S_3$. Similarly, as long as $t_n\in S_3\cup S_4$, the sequence $\{y_n\}=\{x_{n+1}\}$ decreases, until the sequence $\{t_n\}$ steps into $S_1$. Consequently, if a sequence $\{t_n\}$ stays in $S_1$ for all large $n$, we gain a monotonically increasing, bounded sequence $\{x_n\}$, which converges to $B\leq A$. Taking the limit of both sides of $x_{n+1}=a x_n (1-x_{n-1})$, we obtain $B=A$, and consequently $t_0$ is in the region of attraction of $(A,A)$. Similarly, if a sequence gets stuck in $S_3$, it also converges to $(A,A)$. \qed

\bigskip{}

Now we assume $1<a\leq \frac 3 2$ and show that for every $t_0\in S$ the sequence $\{t_n\}$ converges to the nontrivial fixed point $(A,A)$. Combining this fact with the local asymptotic stability of the fixed point (see at the beginning of the following section), Theorem \ref{fotet} is proved for these parameter values.

\begin{all}\label{smalla}
	If $a\in\left(1, \frac {3}{2}\right]$, then $\lim\limits_{n\to \infty} F_a^n(x_0,y_0)=(A,A)$ for every $(x_0,y_0)\in S$.
\end{all}

\noindent \textit{Proof. } It is clear from Proposition \ref{kor}, we only need to consider the case when the sequence $\{t_n\} = \{(x_n, y_n)\}$ goes around the fixed point, not getting stuck in $S_1$ or $S_3$. It means that there exist subsequences $n_k$ and $m_k$, such that
\begin{equation*}
	t_{n_k}, t_{n_k+1}, ..., t_{m_k-2} \in S_1; \quad
	t_{m_k-1} \in S_2; \quad
	t_{m_k}, t_{m_k+1}, ..., t_{n_{k+1}-2} \in S_3; \quad
	t_{n_{k+1}-1} \in S_4
\end{equation*}
for all $k\geq 0$. Clearly $n_k + 2 \leq m_k \leq n_{k+1} -2$ also holds. Without loss of generality, we can assume $t_0 \in S_1$.

Now consider the sequence $\{s_n\}$, where $s_0 = 0$, $s_n = f^n(s_0)$ and $f(x) = a (a-1)(1-x)^2$. Denote by $\{h_n\}$ and $\{g_n\}$ the even and odd indexed subsequences of $\{s_n\}$, i.e. $h_n = s_{2n}$ and $g_n = s_{2n + 1}$. Furthermore introduce the following subsets of $S$:
\begin{equation*}
	\mathcal{H}_k=\{(x,y)\in S : \ h_k\leq x,y \}
	\text{ \ and \ }
	\mathcal{G}_k=\{(x,y)\in S : \ x, y\leq g_k \}.
\end{equation*}
Clearly, $\mathcal{H}_0 = S$.

It is easy to see that, if $t_{n_k} \in \mathcal{H}_k$ then $t_{n_k}, t_{n_k+1}, ..., t_{m_k} \in \mathcal{H}_k$ and because of the inequality
\begin{equation*}
y_{m_k}=ay_{m_k-1}(1-x_{m_k-1})=a^2y_{m_k-2}(1-x_{m_k-2})(1-x_{m_k-1})\leq a^2 A (1-h_k)^2=g_k,
\end{equation*}
$t_{n_k}, t_{n_k+1}, ..., t_{m_k} \in \mathcal{G}_k$ also holds. Similarly, if $t_{m_k} \in \mathcal{G}_k$ then $t_{m_k}, t_{m_k+1}, ..., t_{n_{k+1}} \in \mathcal{G}_k\cap\mathcal{H}_{k+1}$. It follows from the construction, that $\mathcal{H}_{k+1} \subset \mathcal{H}_k$ and $\mathcal{G}_{k+1} \subset \mathcal{G}_k$. Consequently $\{h_n\}$  is increasing and bounded above by $A$, so $\lim_{n\to \infty} h_n=h_{\infty}\leq A$. Similarly $\{g_n\}$ is decreasing and bounded below by $A$, so $\lim_{n\to \infty}g_n=g_{\infty}\geq A$. Therefore, we only need to show that $h_\infty = g_\infty = A$.

It is clear that $h_{\infty}$ and $g_{\infty}$ need to be fixed points of $f_2(x)=f(f(x))$. Observe that $f_2(A)=A$, $f'_2(A)=4(a-1)^2$. Consequently $0<f'_2(A)\leq 1$ for $1< a \leq \frac 3 2$. Furthermore, for $0<x<A$
\begin{equation*}
f''_2(x)=4a^2(1-a)^2 (3f(x)-1) > 0
\end{equation*}
since $f(x)>A\geq \frac 1 3$. We can conclude that $A$ is the only solution of $f_2(x)=x$ in the interval $(0,A]$, so $\lim_{n\to \infty}h_n=A$. From the definition of $g_n$, it is clear, $\lim_{n\to \infty}g_n=A$, too. \qed

\bigskip

In the rest of the paper we assume $a\in \left(\frac{3}{2},2\right]$. For these parameter values, the above argument does not guarantee convergence for every $t \in S$, but we show, it is enough to consider a subset of $S$ later on.

\begin{all} \label{tildeS}
	For every $a\in \left(\frac{3}{2},2\right]$
	the set
	\begin{equation*}
		\tilde S = \left\lbrace (x,y)\in S: \ x,y \in [0.072, 0.8] \right\rbrace
	\end{equation*}
	is invariant, i.e. $F(\tilde S) \subset \tilde S$. Furthermore, for every $t\in S$, there exists $N = N(t)$, such that for every $n>N$, $F^n(t)\in \tilde S$.
\end{all}

\noindent \textit{Proof. } Using the argument of the previous Proposition, we can assume $t_0\in S_1$ and the sequence $\{t_n\} = \{(x_n,y_n)\}$ goes around the fixed point. We need to show that $y_{m_0}\leq \frac{4}{5}$ and $y_{n_1}\geq 0.072$. Since $(x_0,y_0)\in S$ implies $y_1 = a y_0 (1-x_0) \leq a y_0 = a x_1$, we can also assume that $y\leq ax$ for every $t \in S$. 

For $y_{m_0}$ we have to find the maximum of $a^2y(1-y)(1-x)$ assuming $(x,y)\in S_1$ and $y\leq ax$. Since $y(1-y)$ is increasing on $[0,A]$ we are looking for the maximum of $f(x)=a^3x(1-ax)(1-x)$ on $[0,\frac{A}{a}]$ and $g(x)=(a-1)(1-x)$ on $\left[\frac{A}{a},A\right]$. The maxima of $g(x)$ and $f(x)$ are $\frac{3}{4}$ and $\frac{4}{3\sqrt{3}}$, respectively, so $y_{m_0}\leq \frac{4}{5}$. 


Similarly, for every $a\in [1.5,2]$ we are looking for the minimum of $a^2y(1-y)(1-x)$ on $S_3$, assuming $x,y\leq \frac{4}{5}$. It is easy to see that this is $0.072$. \qed

\bigskip{}

We apply this proposition in the computer assisted part of the proof, since it is useful to exclude a small neighbourhood of the trivial fixed point $(0,0)$, as we see it later. For more general results on absorbing sets like $\tilde S$, the reader is referred to \cite{garab2018}.


\section{Attracting neighbourhood with linearisation} \label{seclin}

In this section using the linearisation of map (\ref{falap}), for a fixed parameter $a\in \left(\frac 3 2,2\right)$, we give a neighbourhood $\mathcal{N}(a)$ around $(A,A)$, which is inside the region of attraction of this fixed point, i.e. $\lim_{n\to \infty}F_a^n(x_0,y_0)=(A,A)$ for every $(x_0,y_0)\in \mathcal{N}(a)$.

Introducing the new variables $u=x-A$ and $v=y-A$, map (\ref{falap}) can be written in the following form:
\begin{equation}\label{falap2}
	\left( \begin{array}{c}
		u \\
		v
	\end{array} \right)
	\mapsto
	J(a)
	\left( \begin{array}{c}
		u \\
		v 
	\end{array} \right)
	+
	f_a(u,v),
\end{equation}
where
\begin{equation*}
	J(a)=
	\left( \begin{array}{cc}
	0 & 1  \\
	1-a & 1
	\end{array} \right), \quad
	f_a(u,v)=
	\left(
	\begin{array}{c}
	0 \\
	-a u v 
	\end{array}
	\right)
\end{equation*}

For $a\in \left(\frac 3 2;2\right]$ the eigenvalues of $J(a)$ are $\lambda:=\lambda_{1}(a) = \overline{\lambda_{2}(a)} = \frac{1+i \sqrt{4a-5}}{2}$ and the corresponding eigenvectors are $q_{1,2}(a)=(1, \lambda_{1,2}(a))$. It is easy to see, that $|\lambda_i(a)| < 1$ for $a\in (1,2)$, $i=1,2$, $|\lambda_i(2)| = 1$ and $|\lambda_i(a)| > 1$ for $a>2$, where $i=1,2$. Introduce the notation $q=q(a)=q_1(a)$ and denote by $p=p(a)$ the eigenvector of the transposed matrix $J^T(a)$ corresponding to $\overline{\lambda(a)}$, normalized to $\langle p,q\rangle=1$, where $\langle (a_1,a_2),(b_1,b_2)\rangle=\sum_{i=1}^2 \bar a_i b_i$, $(a_1,a_2),(b_1,b_2)\in\mathbb{C}^2$. We obtain $\bar p=d(\lambda-1, 1)$, where $d=d(a)=(2\lambda(a)-1)^{-1}$.

Introduce the vector $U=(u,v)^T$ and the complex variable $z=\langle p,U\rangle$. The variable $U$ can also be expressed by $z$:

\begin{equation*}
\left( \begin{array}{c}
u \\
v
\end{array} \right)
=U=q(a)z+\overline{q(a)}\bar z=
\left( \begin{array}{c}
z+\bar z \\
\lambda(a)z+\overline{\lambda(a)} \bar z 
\end{array} \right).
\end{equation*}
Moreover, map (\ref{falap2}) can be written in the following form:
\begin{equation*}
z \mapsto \langle p(a),J(a)U+f_a(U)\rangle= \lambda(a)z+d(a)g_a(z),
\end{equation*}
where $g_a(z)=g_a(z,\bar z)=-a(z+\bar z)(\lambda (a) z +\overline{\lambda(a)} \bar z)$ is a real-valued function.

At first we use the map
\begin{equation} \label{zmap}
z \mapsto \lambda(a)z+d(a)g_a(z) :=G(z)
\end{equation}
without further transformation to construct $\mathcal{N}(a)$.


\begin{all} \label{lin}
	For every $a\in \left( \frac{3}{2},2 \right)$ define $\varepsilon(a)$ by
	\begin{equation*}
		\varepsilon(a)=\frac{\sqrt{4a-5}(1-\sqrt{a-1})} {a(2\sqrt{a-1}+1)}\cdot \frac{\sqrt{4a-5}}{\sqrt{a+1}}.
	\end{equation*}
Then the set 
	\begin{equation*}
		\mathcal{N}(a)=\left\lbrace (x,y)\in S:\ |x-A|,|y-A|\leq \varepsilon(a) \right\rbrace
	\end{equation*}
	is in the region of attraction of the fixed point $(A,A)$ of $F_a$.
\end{all}

\noindent\textit{Proof.} At first we show, there exists a $\zeta_0>0$, such  that $|G_a(z)|<|z|$ for every $0<|z|<\zeta_0$. If such a $\zeta_0$ exists, it is clear that the open ball $B_{\zeta_0}^\circ$ around the origin is invariant and we show that every point of $B_{\zeta_0}^\circ$ tends to the origin. Let $z_0$ be an arbitrary point from $B_{\zeta_0}^\circ$ and consider the nonnegative, strictly decreasing sequence $\{|z_n|\}_{n=0}^{\infty}$, where $z_{n+1}=G_a(z_n)$. This sequence can converge only to a fixed point of the continuous map $r\mapsto \max_{|\zeta|=r} |G_a(\zeta)|$, which is, inside $B_{\zeta_0}^\circ$, solely $r=0$.

Estimate the right hand side of the map (\ref{zmap}). Using $|\lambda(a)|=\sqrt{a-1}$, $|d(a)|=\frac 1 {\sqrt{4a-5}}$ and $|g_a(z)|\leq a (2|\lambda(a)|+1)|z|^2=a (2\sqrt{a-1}+1)|z|^2$, we obtain
$$
|\lambda(a)z+d(a)g_a(z)|\leq |z|\left(\sqrt{a-1}+ \frac {a(2\sqrt{a-1}+1)} {\sqrt{4a-5}}|z| \right) <|z|,
$$
for every $z\neq 0$, provided $|z|<\zeta_0:= \frac{\sqrt{4a-5}(1-\sqrt{a-1})} {a(2\sqrt{a-1}+1)}$.

To obtain an estimation of the real variables $u,v$, we use the expression $z=\langle p (a),U\rangle=d((\lambda-1)u+v)$. Supposing $|u|,|v|\leq \delta$, we obtain
$$
|z|\leq |d||(\lambda-1)u+v|=|d| \sqrt{\frac {4a-5} 4 u^2+\left(v-\frac u 2\right)^2}\leq \delta \frac{\sqrt{a+1}}{\sqrt{4a-5}},
$$
therefore, if $\delta\leq \zeta_0 \frac{\sqrt{4a-5}}{\sqrt{a+1}}$, then $|z|\leq \zeta_0$. Set $\varepsilon(a) = \zeta_0 \frac{\sqrt{4a-5}}{\sqrt{a+1}}$. Then points, whose coordinate satisfy $|u|,|v|\leq \varepsilon(a)$, are in the region of attraction of the fixed point. \qed

\bigskip

It is easy to see the set $\mathcal{N}$ shrinks to the fixed point $(A,A)$ as $a$ tends to $2$, since $\lim_{a\to 2}\varepsilon(a)=0$. Consequently, close to the critical parameter value, the neighbourhood, obtained by linearisation is not suitable for reliable numerical methods. In fact, the smaller the neighbourhood, the less efficient, and more time consuming the numerical part of the proof. Furthermore, the linearisation does not provide an attractive neighbourhood at the critical parameter value $a=2$, therefore, we need an other approach to construct a neighbourhood $\mathcal{N}$ for parameter values close to $2$.

In the subsequent section we use the normal form of Neimark--Sacker bifurcation and create a neighbourhood whose size is independent of $a$. Actually, the first method with the linearisation become rather compute-intensive at about the parameter range $(1.99,2)$, but we will use the second technique with the normal form in a bigger parameter range, namely for $a\in (1.95,2]$. The normal form technique provides a significantly larger neighbourhood than the first method can do for parameters close to the critical value, so the second method is more efficient even for $a \in (1.95,1.99]$, too.

\section{Transforming to normal form}\label{secnorm}

In this section, first, we give a general method to construct an attractive neighbourhood around a fixed point, which undergoes a supercritical Neimark--Sacker bifurcation at $a_0$. This neighbourhood is suitable for parameters close to the critical value $a_0$, i.e. for $a\in [a_0-\beta_0,a_0]$ with some fixed $\beta_0>0$. We follow the steps of finding the normal form of the Neimark–Sacker bifurcation, according to Kuznetsov \cite{kuzn}.

Suppose, we have a map
\begin{equation}\label{gen1}
	x \mapsto F_a(x),
\end{equation}
where $x \in \mathbb{R}^2$, $F_a$ is smooth and $a\in \mathbb{R}$ is the parameter. Furthermore, we have a fixed point $\tilde x = \tilde x(a)$, which undergoes a supercritical Neimark--Sacker bifurcation at $a_0$. Fix some $\beta_0>0$.
According to Kuznetsov, if
$|\lambda(a)| < 1$ for all $a\in [a_0-\beta_0,a_0)$, then the map (\ref{gen1}) can be transformed into the following form
\begin{equation} \label{gen2}
z \mapsto G(z) = \lambda(a) z + G_2(z,a),
\end{equation}
where $z\in \mathbb{C}$, and $G_2$ is smooth. (Compare Section \ref{seclin}.)

We can write the smooth $G(z)$ as a formal Taylor series in two complex variables ($z$ and  $\bar z$):
\begin{equation}
G(z)= \lambda(a) z+\sum_{k+l=2}^{4}\frac{g_{kl}}{k!\, l!}z^k \bar z^l+R_1,
\end{equation}
where $g_{kl} = g_{kl}(a)$ and $R_1= R_1(z, \bar z, a)= O(|z|^5)$. Then, with smooth and invertible functions, we transform the map (\ref{gen2}) into the normal form of the bifurcation:
\begin{equation} \label{wmap}
w \mapsto \lambda(a) w + c_1(a)w^2 \bar w + R_2,
\end{equation}
where $R_2= R_2(w, \bar w, a) = O(|w|^4)$.
If we show that there exists $\rho_0>0$, such that for every $0<|w|<\rho_0$ and $a\in [a_0-\beta_0, a_0]$ the following holds
\begin{equation} \label{contr}
|\lambda(a) w + c_1(a)w^2 \bar w + R_2|<|w|,
\end{equation}
then we obtain that $B_{\rho_0}=\{w: |w|<\rho_0\}$ is in the region of attraction of the fixed point $0$ of the map (\ref{wmap}). Since $|\lambda|\leq 1$ for $a\leq a_0$, and the bifurcation is supercritical, i.e. $\Re\frac{{c}_1(a_0)}{\lambda(a_0)}<0$, it is easy to see, that inequality (\ref{contr}) holds for all sufficiently small $\rho_0$ and $\beta_0$.

Our aim is to obtain an explicit value for $\rho_0$ assuming $\beta_0$ is given. Furthermore, $\rho_0$ need to be as big as possible, because of the computer assisted part of the proof. Consequently, the estimation of the higher order terms ($R_2$) is the most essential part of the method, just like in the linearised case. Note that, in the end we need to derive a $\{z: |z|<\varepsilon\}$--type neighbourhood, related to the original map (\ref{gen2}).


To obtain the normal form, we look for a smooth invertible function $h:\mathbb{C} \to \mathbb{C}$ in a neighbourhood of $0 \in \mathbb{C}$ which transforms the map (\ref{gen2}) with the new coordinate $w=h^{-1}(z)$ into the following form:
\begin{equation} \label{highw}
w \mapsto h^{-1}(G(h(w))) = \lambda(a) w + c_1(a)w^2 \bar w + R_2.
\end{equation}
According to Kuznetsov, such a function can be found in the form:
\begin{equation}\label{h}
h(w)=w+\frac{h_{20}}{2}w^2+h_{11}w \bar w+\frac{h_{02}}{2} \bar w^2+\frac{h_{30}}{6}w^3+\frac{h_{12}}{2}w \bar w^2+\frac{h_{03}}{6} \bar w^3,
\end{equation}
where $h_{ij} = h_{ij}(a)$.
To this transformation we need the non-resonance condition
\begin{equation*}
\left( \frac{\lambda(a_0)}{|\lambda(a_0)|} \right)^k \neq 1, \qquad  \textrm{where \ \ } k \in {1,2,3,4}.
\end{equation*}

Clearly $h$ has an inverse in a small neighbourhood of $0 \in \mathbb{C}$, and  $h^{-1}$ can be written in the following form:
\begin{equation}\label{hi}
h^{-1}(z)=h^{-1}_0(z)+R_3,
\end{equation}
where
\begin{equation*}
h^{-1}_0(z) = z+\sum_{2\leq k+l\leq 4}\tilde h_{kl} z^k \bar z^l,
\end{equation*}
$R_3= R_3(z, \bar z, a)= O(|z|^5)$ and $\tilde h_{kl} = \tilde h_{kl}(a)$. The coefficient $\tilde h_{kl}$ can be obtained by substituting $w=h^{-1}(z)$ into $z=h(w)$ and equating the coefficients of the same type up to fourth order. The $h_{ij}$ was obtained in a similar manner: we need to choose the coefficients so that the second and third order terms (apart from $w^2\bar w$) of $h^{-1}(G(h(w)))$ are eliminated.
The formulas can be found in the Appendix. Notice that $h_{kl}$ and consequently $\tilde h_{kl}$ depend only on the at most third order terms of $G$.

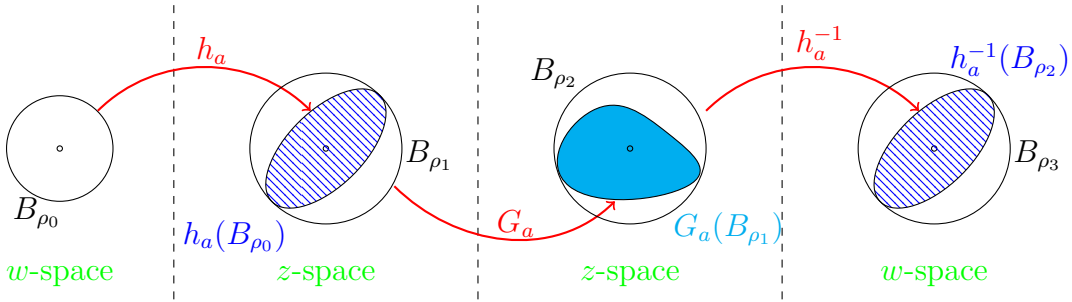
\begin{figure}[h]
	\centering
	\begin{tikzpicture}
	
	\draw[dashed] (3.5,0) -- (3.5,4);
	\draw[dashed] (7.5,0) -- (7.5,4);
	\draw[dashed] (11.5,0) -- (11.5,4);
	
	\draw (2,2) circle[radius=1pt];
	\draw (2,2) circle (.7cm);
	\path (1.7,1.5) node[below] {$B_{\rho_0}$};
	\node[green] at (2,.3) {$w$-space};
	
	\draw [->, thick, red] (2.5,2.5) to[out=45,in=135] (5.3,2.5);
	\node[red] at (4,3.3) {$h_a$};
	
	\draw (5.5,2) circle[radius=1pt];
	\draw (5.5,2) circle (1cm);
	\path (4.3,1.2) node[below, blue] {$h_a(B_{\rho_0})$};
	\draw[shift={(5.5 cm,2 cm)},rotate=45, pattern=north west lines, pattern color=blue] (0,0) ellipse (1cm and .5cm);
	\path (6.4,1.9) node[right] {$B_{\rho_1}$};
	\node[green] at (5.5,.3) {$z$-space};
	
	\draw [->, thick, red] (6.4,1.5) to[out=-45,in=-135] (9.3,1.3);
	\node[red] at (8,1) {$G_a$};
	
	\fill[cyan, draw=black] plot [smooth cycle, tension=1] coordinates {(8.7,1.5) (10.25,1.5) (10,2.2) (9,2.5)};
	\draw (9.5,2) circle[radius=1pt];
	\draw (9.5,2) circle (1cm);
	\path (10.8,1.3) node[below, cyan] {$G_a(B_{\rho_1})$};
	\node[] at (8.5,3) {$B_{\rho_2}$};
	\node[green] at (9.5,.3) {$z$-space};
	
	\draw [->, thick, red] (10.5,2.5) to[out=45,in=135] (13.3,2.5);
	\node[red] at (12,3.4) {$h_a^{-1}$};
	
	\draw (13.5,2) circle[radius=1pt];
	\draw (13.5,2) circle (1cm);
	\path (14.5,2.8) node[above, blue] {$h_a^{-1}(B_{\rho_2})$};
	\draw[shift={(13.5 cm,2 cm)},rotate=45, pattern=north west lines, pattern color=blue] (0,0) ellipse (1cm and .5cm);
	\path (14.4,1.9) node[right] {$B_{\rho_3}$};
	\node[green] at (13.5,.3) {$w$-space};
	
	\end{tikzpicture}
\caption{The size of the domains of $h_a$, $G_a$ and $h_a^{-1}$ if $|w|<\rho_0$ is assumed}\label{r12}
\end{figure}

First we will give a finite-order polynomial estimation on the functions $G_a$, $h_a$ and $h_a^{-1}$:
\begin{eqnarray} 
|h(w)| &\leq & |w|+h_2 |w|^2 + h_3 |w|^3,
\nonumber \\
|G(z)| &\leq  & |z|+g_2 |z|^2+g_3 |z|^3 + g_4 |z|^4 + R_{10} |z|^5,
\nonumber \\
|h^{-1}(z)| &\leq  & |z| + \tilde h_2 |z|^2 + \tilde h_3 |z|^3 + \tilde h_4 |z|^4 + R_{30} |z|^5,
\nonumber
\end{eqnarray}
where the coefficients are independent of $a$. With them we can give an estimation on $R_2$, i.e. the higher order terms of the composition $h^{-1}(G(h(w)))$. Clearly the Taylor expansion of $h$ is finite, but generally the other two Taylor expansions have infinitely many terms. So the at least fifth order terms are estimated in $R_{10}|z|^5$ and $R_{30}|z|^5$. For the lower order terms we have explicit formulae and they could be estimated by interval arithmetic. As for the higher order terms it is essential to be able to say how large can be the moduli of $h(w)$, $G(h(w))$ and $h^{-1}(G(h(w)))$ in (\ref{highw}) if $|w|<\rho_0$ is assumed, since the estimation of the remaining terms of a Taylor expansion highly depends on the size of the neighbourhood on which it need to be valid. The radii $\rho_1$, $\rho_2$ and $\rho_3$ must be chosen so that $h_a(B_{\rho_0})\subset B_{\rho_1}$, $G(B_{\rho_1})\subset B_{\rho_2}$ and $h^{-1}(B_{\rho_2})\subset B_{\rho_3}$ (see figure \ref{r12}), consequently during the study of $G$ and $h^{-1}$ we can assume that the domains are in $B_{\rho_1}$ and $B_{\rho_2}$ respectively. 

After gaining an estimation on $R_2$ we show that inequality (\ref{contr}) holds for $0<|w|\leq\rho_0$. From this result a neighbourhood in the $z$-plane can easily be obtained: the set $B_\varepsilon=\{z:|z|\leq\varepsilon\}$ is inside the attractive neighbourhood of the fixed point of the map (\ref{gen2}) if $h^{-1}(B_{\varepsilon}) \subset B_{\rho_0}$, i.e. $B_\varepsilon$ is mapped inside the region of attraction of the map (\ref{wmap}).

Here, we emphasise that for our calculations the only thing we need to know is the at most fourth order terms of the function $G(z)$ and an $R_{10}|z|^5$--type estimation of the at least fifth order terms of $G(z)$.

\medskip

Until this point in the section we describes our method for a general $F_a(x)$. Now, we turn our attention to the specific $F_a(x,y) = (y, ay(1-x))$ from (\ref{falap}).

The main results of this section are the following two propositions. We prove only Proposition \ref{b1} as the whole argument can be repeated to get an attracting neighbourhood when only $a\in[1.95,2]$ is assumed. The differences appear only in concrete values in the given estimations. Details of Proposition \ref{b0} can be found on our website \cite{DJ}.

\begin{all} \label{b0}
	For all fixed $a\in[1.95,2]$, the set $\{z\in \mathbb{C} : |z|\leq 0.013 \}$ belongs to the basin of attraction of the fixed point $0$ of $G(z)$.
\end{all}

\begin{all} \label{b1}
	For all fixed $a\in[1.995,2]$, the set $\{z\in \mathbb{C} : |z|\leq 0.014 \}$ belongs to the basin of attraction of the fixed point $0$ of $G(z)$.
\end{all}

\noindent
\textit{Proof.} Throughout the proof we suppose $a\in [a_0-\beta_0,a_0]$, where $\beta_0=0.005$ and $a_0=2$. In our calculations we use symbolic computation and built-in interval arithmetic tools of Wolfram Mathematica v. 11.


\subsection{Estimation of the lower order terms in $G_a$, $h_a$ and $h_a^{-1}$}

Throughout this section we need estimation of the coefficients of the lower order terms in $G$, $h$ and $h^{-1}$, such that these estimations are independent of $a\in [1.995,2]$.
We use interval arithmetic tools to compute them for $a\in [1.995,2]$.

In our particular case the function $G(z)$ can be written in the following form:
$$
G(z)= \lambda(a) z+\sum_{k+l=2}\frac{g_{kl}}{k!\, l!}z^k \bar z^l,
$$
since $G$ has only at most second order terms. Furthermore, we use (\ref{h}) and (\ref{hi}).
We look for constants satisfying the following inequalities:
\begin{align*}
g_2 & \geq  \max_{a \in [a_0-\beta_0, a_0]} \left( \frac{|g_{20}|}{2}+|g_{11}|+\frac{|g_{02}|}{2} \right),  \\
h_2 & \geq  \max_{a \in [a_0-\beta_0, a_0]} \left( \frac{|h_{20}|}{2}+|h_{11}|+\frac{|h_{02}|}{2} \right), \\
h_3 & \geq   \max_{a \in [a_0-\beta_0, a_0]} \left( \frac{|h_{30}|}{6}+\frac{|h_{12}|}{2}+\frac{|h_{03}|}{6} \right),  \\
\tilde h_n & \geq \max_{a \in [a_0-\beta_0, a_0]} \left(\sum_{i+j=n} |\tilde h_{ij}| \right), 
\end{align*}
where $n=2,3,4$.
With interval arithmetic it can be shown that
$g_2=3.47$, $h_2=2.9$, $h_3=4.7$, $\tilde h_2=h_2=2.9$, $\tilde h_3=8.2$ and $\tilde h_4=30$
fulfil the requirements. From the definition of these constants we obtain the following estimations:
\begin{align} \label{max0}
|h(w)|&  \leq h^{max}(|w|)  :=  |w|+h_2 |w|^2 + h_3 |w|^3, \nonumber \\
|G(z)| & \leq G^{max}(|z|)  :=  |z|+g_2 |z|^2, \\
|h^{-1}_0(z)| & \leq \tilde h^{max}_0(|z|)  :=  |z|+\tilde h_2 |z|^2 + \tilde h_3 |z|^3 + \tilde h_4 |z|^4, \nonumber \\
|h^{-1}(z)| & \leq \tilde h^{max}(|z|)  :=  |z| + \tilde h_2 |z|^2 + \tilde h_3 |z|^3 + \tilde h_4 |z|^4 + R_{30} |z|^5, \nonumber
\end{align}
if in the last equation $R_{30}$ satisfies $|R_3|\leq R_{30}|z|^5$. We will determine $R_{30}$ later.

From the definition of $h_2$ and $h_3$ we also get:
\begin{equation} \label{hiz}
|w|- h_2 |w|^2 - h_3 |w|^3 \leq |h(w)|,
\end{equation}
Consequently, assuming $|w|\leq \rho$, we can make the following estimation:
\begin{equation} \label{wkz}
|w| \leq \eta(\rho)|h(w)|,
\end{equation}
with
\begin{equation*}
	\eta(\rho) := \frac 1 {1 - h_2 \rho - h_3 \rho^2}.
\end{equation*}


We choose $\rho_0=0.015$, $\rho_1=h^{max}(\rho_0)$, $\rho_2=G^{max}(\rho_1)$ and from (\ref{hiz}) it is clear that $B_{\rho_2}$ can not be mapped outside of the circle with radius $0.018$, consequently this value is a suitable choice for $\rho_3$.



\subsection{The domain of $h$ and $h^{-1}$}

Now, we show that $h$ is injective in $\overline {B_{1/9}} \subset \mathbb{C}$, and $h^{-1}$ is defined on $\overline {B_{1/16}}$. Let $z \in \mathbb{C}$, $a \in [1.995, 2]$ be fixed, and denote $ H_{a,z}: \mathbb{C} \ni w \mapsto w+z-h(w) \in  \mathbb{C} $. With this notation $H_{a,z}=w$ if and only if $h(w)=z$.
\begin{align*} 
|H_{a,z}(w_1)-H_{a,z}(w_2)| & =  |w_1-h(w_1)-w_2+h(w_2)|    \\
& \leq |w_1-w_2| \left( h_2 (|w_1|+|w_2|)+h_3  (|w_1|^2+|w_1|\cdot|w_2|+|w_2|^2) \right)   
\end{align*}
If $|w|\leq \delta_1$ and $|z| \leq \delta_2$, then 
\begin{equation*}
|H_{a,z}(w_1)-H_{a,z}(w_2)| \leq |w_1-w_2|  (2 \delta_1 h_2+3 \delta_1^2 h_3)
\end{equation*}
and
\begin{equation*}
|H_{a,z}(w)| \leq \delta_2 + \delta_1^2 h_2 + \delta_1^3 h_3.
\end{equation*}

Choosing $\delta_1=\frac{1}{9}$ and $\delta_2=\frac{1}{16}$ the map $H_{a,z}$ is a contraction mapping $\overline {B_{1/9}}$ into itself. Consequently for every $z \in \overline {B_{1/16}}$ there exists only one $w=w(z) \in \overline {B_{1/9}}$ such that $h(w(z))=z$, i.e. $h^{-1}$ is defined on $\overline {B_{1/16}}$.

It is clear, that $\rho_0,\rho_3<\delta_1$ and $\rho_1,\rho_2<\delta_2$, where $\rho_0,\rho_1,\rho_2,\rho_3$ were chosen at the end of the previous subsection.

\subsection{The estimation of the higher order terms in $h^{-1}$}

Now, we turn our attention to the estimation of $R_3$ in (\ref{hi}), which consists of the fifth and higher order terms of $h^{-1}$. We need an estimation $|R_3(z)| < R_{30} |z|^5$, assuming $|z|\leq \rho_2$. But first, we give an estimation of type $|R_3(h(w))| < C |w|^5$, assuming $|w|\leq \rho_3$ (see figure \ref{r12}). Using the definition of $h^{-1}_0$ and $h$, it follows from (\ref{hi}) that
\begin{equation*}
R_3(h(w))=w - h^{-1}_0(h(w))=\sum_{5\leq k+l \leq 12}r_3^{kl}(a) w^k \bar w^l,
\end{equation*}
since it is a twelfth order polynomial of $w$ and $\bar w$, which has only fifth and higher order terms.

Consider the composition
\begin{equation*}
	\tilde R_3(|w|) = \tilde h_0^{max}(h^{max}(|w|)) = \sum_{j=5}^{12} r_3^j |w|^j
\end{equation*}
of the real functions $h^{max}, \tilde h_0^{max}$. It is clear, that $\sum_{k+l=j} |r_3^{kl}| \leq r_3^j$ holds for $5 \leq j \leq 12$.
Consequently,
\begin{equation*}
\left|R_3(h(w))\right|\leq\sum_{5\leq k+l \leq 12}\left|r_3^{kl} (a)\right| |w|^{k+l} \leq \sum_{5\leq j \leq 12}r_3^{j} |w|^{j} \leq \sum_{5\leq j \leq 12}r_3^{j} \rho_3^{j-5} |w|^5,
\end{equation*}
assuming $|w|\leq \rho_3$.
Using (\ref{wkz}), we gain $|w|<\eta(\rho_3)|z|$, and
\begin{equation*}
|R_3(z)| \leq \sum_{5\leq j \leq 12}r_3^{j} \rho_3^{j-5} (\eta(\rho_3))^5|z|^5 \leq 1070|z|^5,
\end{equation*}
therefore, $R_{30} = 1070$ is a suitable choice. 



\subsection{The estimation of the higher order terms in the normal form}

Now we turn our attention to $R_2$, which estimates the at least fourth order terms of $h^{-1}(G(h(w)))$. To obtain a better estimation, we handle the fourth order terms ($R_{24}$) and the higher order ones ($R_{25}$) separately. Set $R_2(w) = R_{24}(w)+R_{25}(w)$.

The fourth order coefficients $r_2^{kl}(a)$ (where $k+l=4$) of $h^{-1}(G(h(w)))$ can be calculated explicitly; the formulae of $\left|r_2^{kl}(a)\right|$ can be found in the Appendix. With interval arithmetic it can be shown, that $\sum_{k+l=4}\left|r_2^{kl}(a)\right| \leq 40$, consequently $\left|R_{24}(w)\right| \leq 40 |w|^4$.

As for the higher order terms, we use $h^{max}, G^{max}$ and $\tilde h^{max}$, similarly to the estimation of $R_3$. Consider the composition
\begin{equation*}
	\tilde R_2\left(|w|\right) = \tilde h^{max}(G^{max}(h^{max}(|w|))) = \sum_{j=1}^{30} r_2^j |w|^j.
\end{equation*} 
It is clear, that for $|w| \leq \rho_0$ the following holds: 
\begin{equation*}
	|R_{25}(w)| \leq \sum_{k=5}^{30} r_2^k(|w|) \rho_0^{k-4} |w|^4 \leq 90 |w|^4.
\end{equation*}

Combining these two results, for $|w|\leq \rho_0$, we obtain
\begin{equation*}
|R_2(w)| \leq |R_{24}(w)|+|R_{25}(w)|\leq 130 |w|^4,
\end{equation*}
and consequently $R_{20}=130$.


\subsection{The attracting neighbourhood}

Now, with our previous estimation on $R_2$ we can finish our proof. Since
\begin{equation*}
|\lambda(a) w + c_1(a)w^2 \bar w + R_2|
\leq  |w| \left( \left| |\lambda | + \tilde c_1 |w|^2 \right| + R_{20}|w|^3 \right),
\end{equation*}
where $\tilde c_1=\frac{|\lambda(a)|} {\lambda(a)}c_1(a)$ and $\lambda=\lambda(a)$, we only need to prove
\begin{equation} \label{kontr}
\left| |\lambda | + \tilde c_1 |w|^2 \right| + R_{20}|w|^3 < 1
\end{equation}
for every $|w|\leq \rho_0$ and $a\in [a_0-\beta_0,a_0]$.

To this end, we show that the following inequality holds with a suitable $R_{4}>0$:
\begin{equation*} \label{rec}
\big| |\lambda| + \tilde c_1 |w|^2 \big| \leq |\lambda| + (\Re\tilde c_1) |w|^2 + R_{4} |w|^3,
\end{equation*}
or equivalently
\begin{equation} \label{rec2}
0 \leq 2 R_4 |\lambda|-(\Im \tilde c_1)^2|w|+2 R_4(\Re \tilde c_1) |w|^2+ R_4^2 |w|^3,
\end{equation}
for every $|w|\leq \rho_0$ and $a\in [a_0-\beta_0,a_0]$. For $a\in [a_0-\beta_0,a_0]$ we can make the following estimations with interval arithmetic on the coefficients in (\ref{rec2}) depending on $a$: $\Re \tilde c_1$ and $\Im\tilde c_1$ are negative, $|\Re\tilde c_1|\leq 2.1$, $|\Im \tilde c_1|\leq 3.5$ and $|\lambda|\geq 0.99$. From this it is clear that for $|w|\leq \rho_0$ the choice $R_4=0.1$ will be suitable. Therefore the left hand side of the inequality (\ref{kontr}) can be written in the following form:
\begin{eqnarray}
\big| |\lambda| + \tilde c_1 |w|^2 \big| + R_{20}|w|^3 & \leq & (|\lambda| + \Re\tilde c_1 |w|^2) + (R_4 + R_{20}) |w|^3\leq
\nonumber \\
& \leq & 1+\left(\Re\tilde c_1+ (R_4 + R_{20}) |w|\right) |w|^2,
\nonumber
\end{eqnarray}
which is less than $1$, provided
$$
|w|< \frac{-\Re \tilde c_1}{R_4 + R_{20}}.
$$
Using the fact that $|\Re\tilde c_1|\geq 2$ we obtain $\frac{-\Re \tilde c_1}{R_4 + R_{20}} > \rho_0$, therefore inequality (\ref{kontr}) holds for every $|w|<\rho_0$. From $|h_{a}^{-1}(z)|<h^{max}_{inv}(|z|)$, inequality $|z|<\varepsilon_G:=0.014$ implies $|w|=|h^{-1}(z)|<\rho_0$, so the proof of Proposition \ref{b1} is complete.
\qed

\bigskip



To obtain a neighbourhood in the real coordinate system $(u,v)$ we use $z=\langle p (a),U\rangle$, just like in the linearised case. 
Note that, the set $\{z\in \mathbb{C} : |z|\leq \varepsilon_G \}$ will be transformed into an ellipse-shaped neighbourhood in the $uv$-plane.

\newpage

\section{Graph representation}\label{secpc}

In the computer assisted part, we follow the method in \cite{bgk}. In this section (for $\frac{3}{2} < a \leq 2$) we associate the map (\ref{falap}) with a directed graph, which reflects the behaviour of the map up to a given resolution. Therefore we can derive properties of our dynamical system through the study of this graph. More precisely our aim is to show with the help of this graph, that every point of $S\setminus \mathcal{N}$ enters into the attracting neighbourhood $\mathcal{N}$ of the nontrivial fixed point constructed in the previous sections.

	Let $D$ be a subset of $\mathbb{R}^n$. A set $\mathfrak{S}$ is called a \textit{cover} of $D$, if the elements of $\mathfrak{S}$ are subsets of $\mathbb{R}^n$ and $\cup_{s\in \mathfrak{S}}s\supset D$.
	Let a map $f: D_f\subset\mathbb{R}^n\to \mathbb{R}^n$, a subset $D\subset D_f$ and a cover $\mathfrak{S}$ of $D$ be given. The directed graph $G(V,E)$ is called a \textit{graph representation} of $f$ on $D$ with respect to $\mathfrak{S}$, if there exists a bijection $i:V\to \mathfrak{S}$, such that the following implication is true for all $u,v\in V$:
	$$ f(i(u)\cap D)\cap i(v)\cap D\neq \emptyset \Rightarrow (u,v)\in E.$$

The meaning of the implication in the previous definition is the following. If we can get with map $f$ from an element $s_1$ of the cover to an other (possibly the same) element $s_2$ of it, i.e. there exists $x\in s_1$ and $y\in s_2$ such that $f(x)=y$, then there is an edge between the vertices corresponding to the two sets, more precisely $(u,v)\in E$ for $s_1=i(u)$ and $s_2=i(v)$. The reverse implication is not necessarily true, namely if there is a directed edge between the vertices $u$ and $v$, it is not sure there exists $x\in s_1$ such that $f(x)\in s_2$, where $s_1$ and $s_2$ are the corresponding sets to $u$ and $v$.

It is easy to see the implication above can be reformulated as follows. For every $u\in V$
\begin{equation}
\label{impl} f(i(u)\cap D)\subset \bigcup_{v\in K_u} i(v)\cap D,
\end{equation}
where $K_u$ denotes the set of vertices, into which there is an edge from $u$ in graph $G$. So the sets corresponding to vertices in $K_u$ need to form a cover of the image of $i(u)$. From this it can be seen the graph representation can be regarded as some kind of upper estimation of the original map $f$. The finer the cover is, the better the graph representation approximates the map. Therefore, if we would like to determine the possible location of the image of a point $P\in \mathbb{R}^n$ under $f$, we can do it with the help of the graph, since $f(P)\in \cup_{v\in K_u} i(v)\cap D$ for $P\in i(u)$. This means iterating $f$ the point $P$ can move forward only along the edges, i.e. it can move from an element of the cover to an other one only if there is an edge between the two vertices corresponding to them. Consequently we can draw conclusions regarded the possible future location of a point studying only the graph. In the following we take the liberty to handle the elements of the cover as vertices and vice versa, omitting the use of $i$.

The construction of the graph representation in our case is the following. For a fixed $k\in\mathbb{N}$ we divide the unit square $[0,1]\times[0,1]$ parallel to the sides into $2^k \times 2^k$ pieces of small closed squares with side length $r=2^{-k}$. According to Proposition \ref{tildeS} we only need to consider the squares lying in $\left[ \frac{1}{16},\frac{7}{8}\right] ^2$. The cover $\mathfrak{S}$ of $\tilde S$ consists of these sets. The small squares correspond to the vertices of the graph. As for the edges, for every small square $s$ we construct a rectangle with reliable numerical methods  which contains $f(s)$. If the rectangle intersects the small square $s_2$, then there is an edge from $s$ to $s_2$. It is clear, that this construction satisfies relation (\ref{impl}). Note that we considered only that part of the rectangle obtained by the numerical method, which lies inside the square $\left[ \frac{1}{16},\frac{7}{8}\right] ^2$, but this is not a restriction, since the studied set $\tilde S$ is invariant under the map (\ref{falap}), so getting out of the unit square is only the consequence of the numerical method and the 'upper estimation' nature of the graph representation. Note also that, instead of map (\ref{falap}) we use the second iterate of it, since the formula is still compact enough not to cause big overestimation in interval arithmetic and it considerably speeds up the calculations. 

In this paper we suppose a graph is always finite. A graph is \textit{strongly connected} if there are $uv$ and $vu$ (directed) paths for every $u\neq v$ vertices of the the graph. We use the following decomposition of a directed graph (see \cite{bangdigraphs}).

\begin{all}\label{sccprop}
	The vertices of a directed graph can be classified and the classes can be ordered such that
	\begin{itemize}
		\item the subgraphs spanned by the classes are strongly connected, and
		\item for every directed edge between these classes, the class of the tail of this edge precedes the class of the head of it,
	\end{itemize}
	moreover the partition above is unique.
\end{all}

The aforementioned classes are called the \textit{strongly connected components} (\textit{SCC}) of the graph. A strongly connected component is called \textit{non-essential}, if it consists of one vertex without loop. Otherwise we call it \textit{essential}.

From the graph representation and from Proposition \ref{sccprop} it is clear what happens to an arbitrary point of $S$ during the iterations. Starting from a small square containing this point it moves to an other (possibly the same) small square along a directed edge. If we are not in an essential SCC we step out of this small square not returning to here afterwards because of the ordering of the SCCs. If we are in an essential SCC it can happen, that the point stays here forever, or the point steps out of this SCC, but in this case it can not return to this SCC any more.

Since during the partition we obtain finitely many small squares and consequently the graph is finite, it is straightforward that for every point of $S$ there exists an essential strongly connected component, which the point enter during the iteration and never leaves it. So it is true for every $x\in S$ that it enters an essential SCC with finitely many steps and stays here afterwards, therefore we only need to study the essential SCCs.

Our aim is to show that those essential SCCs, in which the points of $S$ can get stuck, are in the attracting neighbourhood $\mathcal{N}$ of the fixed point $(A,A)$, which neighbourhood was constructed analytically in the previous sections. It is important to note that, it is possible for some essential SCC that none of the points of $S$ can get stuck here. Actually, this would be the case close to the trivial fixed point $(0,0)$, since it is a saddle; that shows the necessity of Proposition \ref{tildeS} and $\tilde S$.



As a next step we refine the partition as follows. We divide the small squares into four smaller squares, that have a side length half as long as before, determine their images with reliable numerical methods and construct the SCCs again. Because of the properties of interval arithmetic (inclusion isotonicity: $I_1\subset I_2 \Rightarrow F(I_1)\subset F(I_2)$, where $F$ is the interval-extension of $f$, see \cite{tuc}), if there is an edge between two new small squares, then there must be an edge between their predecessors with the same orientation. We come to the conclusion that during the refinement, an essential SCC can arise only from a former essential SCC, therefore it is really enough to trace merely the essential SCCs. Note that, with the refinements the graph representation 0ecomes a more and more accurate approximation of the represented map, so an essential SCC can fall apart into smaller pieces, and it even can happen that none of the small squares born from a former essential SCC compose a new essential SCC, i.e. this cycle in the graph is only the consequence of the 'upper estimation' nature of the graph representation.
We continue these refinement steps, until all the remaining SCCs are inside the region of attraction of the fixed point obtained in the previous section. If it occurs in finitely many steps our main theorem is proven.

Finally, instead of checking after every refinement, whether the remaining SCCs are in the analytically constructed attracting neighbourhood $\mathcal{N}$, we can remove all the small squares lying entirely in $\mathcal{N}$ before the first refinement. In that case for a fixed $a$, the main theorem will be proved, if the set of the new SCCs will be empty after a refinement. We show the correctness of this method.
\begin{itemize}
	\item If we erase a vertex which is a non-essential SCC, it has no effect at all compared to our former method (when checking after every refinement).
	\item If we remove a whole essential SCC, it also has no substantial effect, because during the checking it always would be in the attracting neighbourhood.
	\item The only significant change happens, when we erase only a proper subgraph of an essential SCC. Consider such an SCC and colour blue the vertices we want to remove (and do not remove them yet). Delete the directed edges stemming from a blue vertex, then form the SCCs (referred to as new SCCs later on) of the new graph and order them such that the blue vertices are at the end of the ordering. (It can be done, since there are no edges from coloured to uncoloured vertices.)
	\begin{itemize}
		\item[--] An uncoloured vertex can be in a new essential SCC; in that case they remain under study after the removal of the blue vertices, just as they would be in the original method. 
		\item[--] However, if an uncoloured vertex is a non-essential SCC it will be erased (as we keep only the essential SSCs), unlike in the method without deleting the vertices of the attracting neighbourhood, but this is not a problem because every point of this vertex enters a new SCC or a blue vertex (because of the ordering) in finitely many steps, so this vertex really can be deleted.
	\end{itemize}
	
\end{itemize}

Note that the aforementioned method can be regarded as a proof, since the graph problems are finite, so the computer can work on them punctually, moreover the method used during the construction of edges was executed with reliable numerical methods, therefore if we have sufficiently much time, then we could reconstruct by hands the parts which were executed by the computer, and we would come to the same conclusion, if our estimation is as good as the computer's.

\noindent \hrulefill

\noindent \textbf{Algorithm} Proving the global stability of $(A,A)$ for the logistic map

\vspace{-3mm}
\noindent \hrulefill

\noindent \hspace*{.08cm} 1: \textbf{procedure} Log2d

\noindent \hspace*{.08cm} 2: \hspace*{.5cm} $V\leftarrow$ the initial partition \hfill $r=2^{-10}$

\noindent \hspace*{.08cm} 3: \hspace*{.5cm} $E\leftarrow$ the edges \hfill \hfill construct them with reliable num. method

\noindent \hspace*{.08cm} 4: \hspace*{.5cm} $C\leftarrow$ SCC of directed graph($V,E$)

\noindent \hspace*{.08cm} 5: \hspace*{.5cm} 
\textbf{remove} the nonessential SCCs from $V$

\noindent \hspace*{.08cm} 6: \hspace*{.5cm} 
\textbf{remove} the SCC at the origin from $V$ \textbf{if} possible 

\noindent \hspace*{.08cm} 7: \hspace*{.5cm} 
\textbf{remove} the initial attracting neighbourhood from $V$

\noindent \hspace*{.08cm} 8: \hspace*{.5cm} \textbf{repeat}

\noindent \hspace*{.08cm} 9: \hspace*{1cm} $V\leftarrow$ refine($V$) \hfill $r\leftarrow r/2$

\noindent 10: \hspace*{1cm} $E\leftarrow$ the edges

\noindent 11: \hspace*{1cm} $C\leftarrow$ SCC of directed graph($V,E$)

\noindent 12: \hspace*{1cm}
\textbf{remove} the nonessential SCCs from $V$


\noindent 13: \hspace*{.5cm} \textbf{until} $|B_1|=\emptyset$

\noindent 14: \textbf{end procedure}

\noindent \hrulefill

The program code, and the outputs can be found on link \cite{DJ}.

\section{Completion of the proof}

In the previous sections we obtained an attracting neighbourhood and then a method to prove the global stability of the nontrivial fixed point for a fixed $a\in [1.5,2]$.
In this section we show, how to modify our method to handle not only a single value of $[1.5,2]$ but a small subinterval of that, instead. 

Let $[a]=[a^-,a^+]\subset [1.5, 2]$ be a fixed small interval. First, we need a new attracting neighbourhood $\mathcal{N}([a])$, such that for every $a\in [a]$, the attracting neighbourhood $\mathcal{N}(a)$ contains this set, i.e. $\cap_{a\in [a]} \mathcal{N}(a) \supset \mathcal{N}([a])$. To this end, we need to take into consideration the displacement of the fixed point and the change in the size of the neighbourhood. Secondly, during the construction of the edges of the graph representation the number $a$ have to be replaced by the interval $[a]$, since (\ref{falap}) and consequently its second iterate depends on $a$. So while studying the image of a small square $k_1$, we need to study it for every $a\in [a]$, i.e. we take a set of small squares during the estimation of the image set such that they cover $f_a^2(k_1)$ for every $a\in [a]$.


We divide the interval $[1.5,2]$ into subintervals $\mathcal{I}_1 = [1.5,1.95]$, $\mathcal{I}_2 = [1.95,1.995]$ and $\mathcal{I}_3 = [1.995,2]$, then divide further these intervals into smaller subintervals with length $2^{-10}, 2^{-13}$, and $2^{-16}$ respectively (see Table \ref{inttable}).



\begin{table}
	\begin{center}
		\begin{tabular}{|c|c|c|c|c|}
			\hline 
			& parameter & size of slices & shape of $\mathcal{N}$ &    {parameters of $\mathcal{N}$}   \\ 
			\hline \hline
			$\mathcal{I}_1$ & $[1.5,1.95]$ & $2^{-10}$ & rectangle &  $7\cdot 10^{-3}$ \\
			\hline
			$\mathcal{I}_2$ & $[1.95,1.995]$ & $2^{-13}$ & ellipse &   $\varepsilon_G = 0.0138 $\\
			\hline
			$\mathcal{I}_3$ & $[1.995,2]$ & $2^{-16}$ & ellipse &   $\varepsilon_G = 0.0146 $ \\   
			\hline 
		\end{tabular} 
		\caption{The partition of the parameter range}
		\label{inttable}
	\end{center}
\end{table}

For small intervals in $\mathcal{I}_1$ we use the linearised map and the square-shaped neighbourhoods with side length $2\varepsilon(a)$ (Proposition \ref{lin}). It is easy to see the size of this set and the location of the fixed point also changes as $a$ changes. However, it can be shown that $\varepsilon(a)\geq 0.007 $ for every $a\in [1,1,95]$, so considering this value fixed, we need to handle only the displacement of the fixed point. 


For small intervals in $\mathcal{I}_2$ and $\mathcal{I}_3$ we use the bifurcation normal form, therefore, the size of the ellipse-shaped neighbourhood is fixed (Propositions \ref{b0} and \ref{b1}), so we only need to consider the displacement of the fixed point.


\section{The algorithm}

During the calculation of edges of the graph representation we use the second iterate of the original map (\ref{falap}):
\begin{equation}  \label{f2}
\left( \begin{array}{ccc}
x_0 \\
y_0
\end{array} \right)
\mapsto 
\left( \begin{array}{ccc}
x_2 \\
y_2
\end{array} \right)
=
\left( \begin{array}{ccc}
a y_0 (1-x_0)\\
a^2 y_0 (1-x_0)(1-y_0)
\end{array} \right).
\end{equation}
Regarding the examined parameter domain $[a^-,a^+]$ and the sides $[x_i^-,x_i^+]$ and $[y_i^-,y_i^+]$ of the squares as intervals, simply, we could use interval arithmetic tools, such as IntLab to compute the image of a small square. However, the map is quite simple, so we can accelerate this method as follows. Notice that $x_2^-=a^- y_0^- (1-x_0^+)$, so we only need to force the computer to use a downward rounding in order to guarantee that the obtained below estimation is really not larger than the possible first coordinates of the image of any point from the initial square. Similarly we can estimate $x_2^+$, but at this time we use upward rounding. As for the $y_2^-$ and $y_2^+$, remark that instead of $y_0^-(1-y_0^+)$ in $y_2^-$ we can use $\min\{ {y_0^-(1-y_0^-),y_0^+(1-y_0^+)}\}$ because $y_0$ denotes the same number in expression (\ref{f2}), and the function $x(1-x)$ is monotone on intervals which do not contain $\frac 1 2$ in the interior (and it is fulfilled in the partition). In $y_2^+$ we replace the minimum by maximum, after that we proceed just like in the case of $x_2$.

We implemented our program in MATLAB, and used the built-in \textit{digrap} function to construct the directed graph from the edge list and the \textit{conncomp} function to divide the graph into strongly connected components.

\begin{figure}[h]
	\centering
	\subfigure[]{
		\includegraphics[width=5cm]{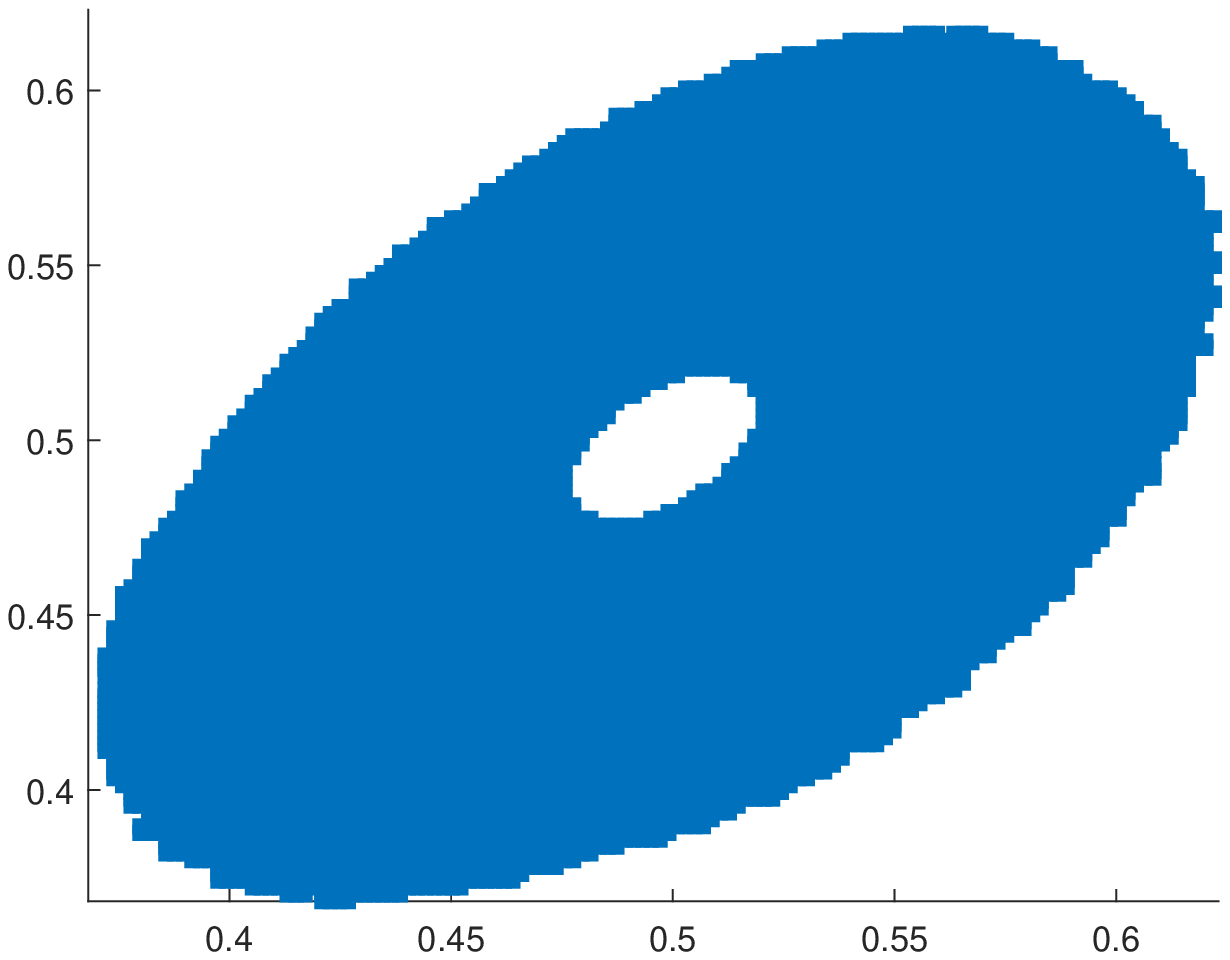}
		\label{fig:subfig1}
	}
	\subfigure[]{
		\includegraphics[width=5cm]{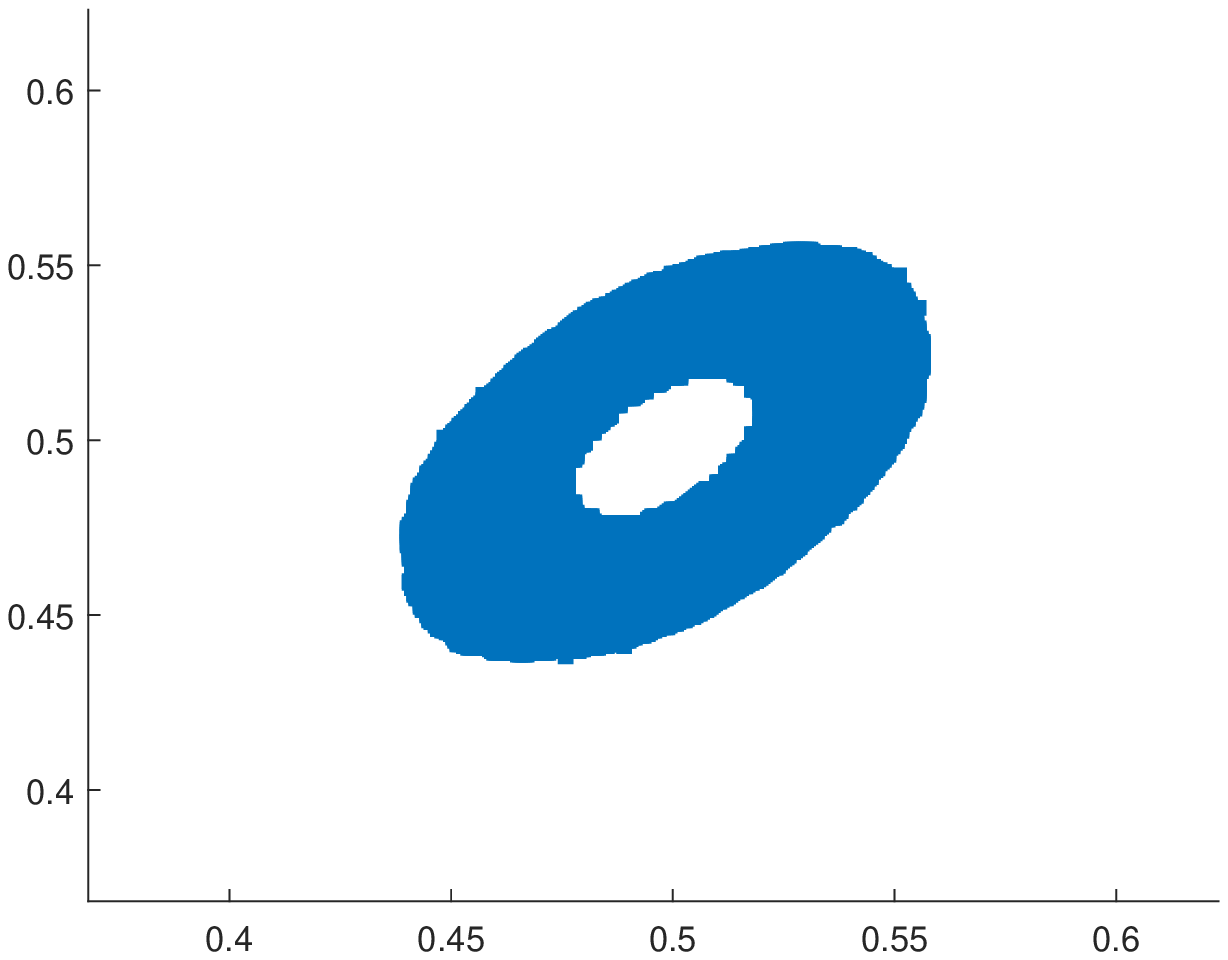}
		\label{fig:subfig2}
	}
\subfigure[]{
	\includegraphics[width=5cm]{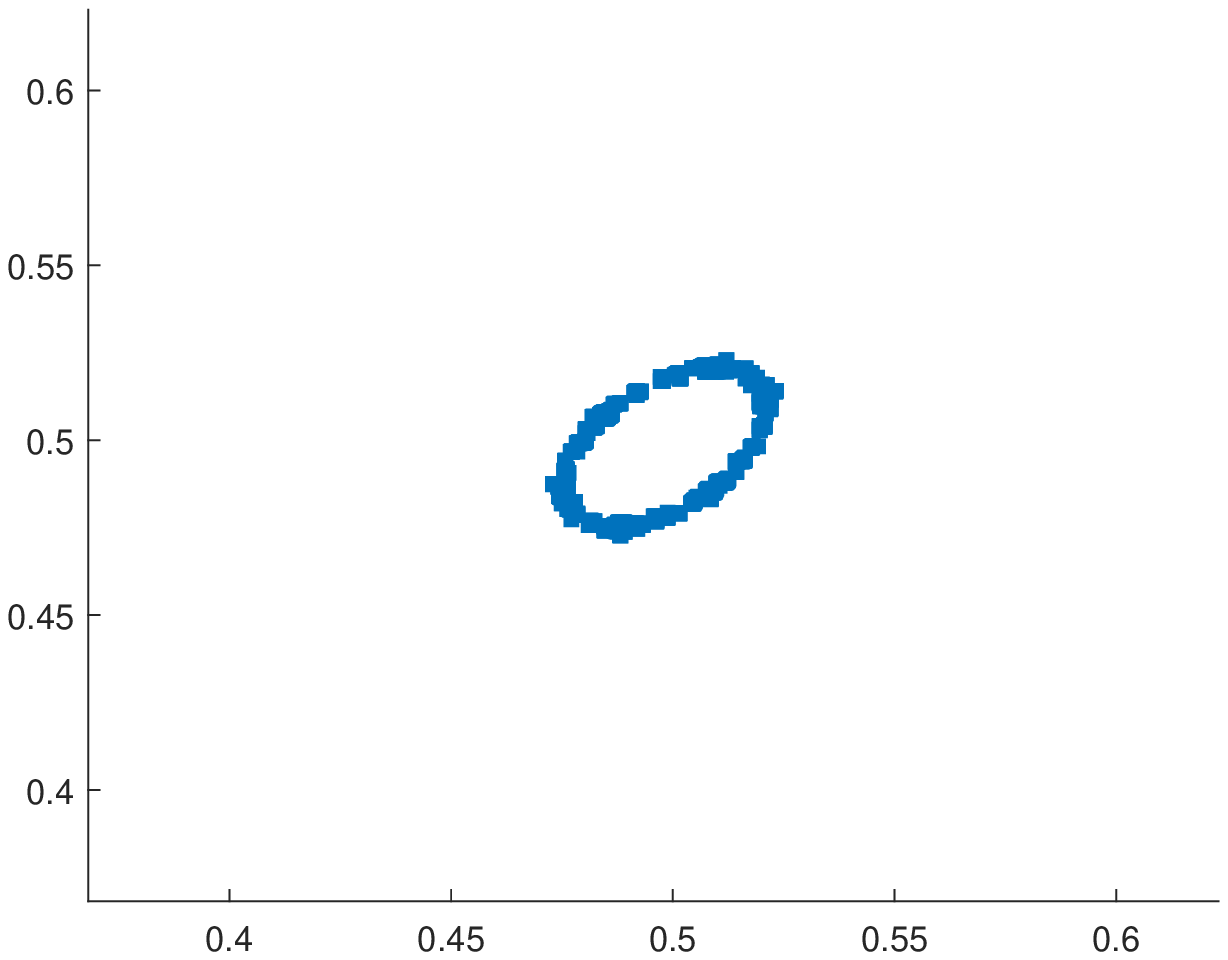}
	\label{fig:subfig3}
}
	\caption{For $[a_-,a_+]=[2-\frac{61}{2^{13}}, 2-\frac{60}{2^{13}}]$ the remaining vertices before the first refinement \subref{fig:subfig1}, after $2$ refinements \subref{fig:subfig2} and after $4$ refinements \subref{fig:subfig3}}
	\label{ex}
\end{figure}

Now, we can run our algorithm with parameters summarised in Table \ref{inttable}. As an example, for the parameter slice $[2-\frac{61}{2^{13}}, 2-\frac{60}{2^{13}}]$ 
we show the evolution of the remaining SCCs during the first $4$ iterations on Figure \ref{ex}.


The program ran successfully, therefore we established the nontrivial fixed point is globally attracting for $a\in[1.5,2]$. Combining this with Proposition \ref{smalla} and the asymptotic stability for $a\in (1,2]$ the proof of Theorem \ref{fotet}
is completed.

\section{Acknowledgement}

I would like to sincerely thank Professor Tibor Krisztin for his useful suggestions. This research was supported by the Hungarian Scientific Research Fund, Grant No. K 109782.

\newpage
\section{Appendix}

\begin{equation*}
	\lambda(a) = \frac{1}{2}+\frac{1}{2}i \sqrt{4a-5}
\end{equation*}

$$
g_{20}(a)=-a + \frac{i}{\sqrt{4a-5}}, \qquad g_{11}(a)= \frac{i a}{\sqrt{4a-5}}, \qquad g_{02}(a)=a + \frac{i}{\sqrt{4a-5}}
$$

$$
h_{20}(a) = \frac{4a}{4a-5 + i \sqrt{4a-5}},
\qquad
h_{11}(a) = \frac{4ia}{\sqrt{4a-5} \left(- i + \sqrt{4a-5} \right)^2},
$$

$$
h_{02}(a)=\frac{a \left( i-\frac{1}{\sqrt{4a-5}} \right)}{i-ia + \sqrt{4a-5}},
$$

\begin{equation*}
h_{30}(a) = -\frac{12i a^2}{-2\left( 5i+\sqrt{4a-5}\right) + a \left( 13i - \sqrt{4a-5} + 2a \left( -2i +\sqrt{4a-5}\right)\right)},
\end{equation*}

\begin{equation*}
h_{12}(a) =  \frac{16a^2 \left( 2-2i\sqrt{4a-5} + a\left( -5 + 2a -\sqrt{4a-5} \right) \right)}{\sqrt{4a-5}\left( -i + \sqrt{4a-5} \right)^4 \left( -7i - 3\sqrt{4a-5} + a \left( 7i - ia + 2\sqrt{4a-5} \right) \right)}
\end{equation*}

\begin{equation*}
h_{03}(a) =  \frac{96 a^2 \left( -2 + a + i\sqrt{4a-5} \right)}{\sqrt{4a-5} \left( 1 + i\sqrt{4a-5} \right)\left( a - 1  + i\sqrt{4a-5} \right) \left( 12 i (a-1) - 7\sqrt{4a-5} + \sqrt{(4a-5)^3} \right)}
\end{equation*}

$$
h_{inv}^{20}(a)=-\frac{h_{20}}{2}, \qquad h_{inv}^{11}(a)=-h_{11}, \qquad h_{inv}^{02}(a)=-\frac{h_{02}}{2}
$$

$$
h_{inv}^{30}(a)=\frac{1}{6}(3 h_{20} ^2-h_{30} +3 h_{11} \bar h_{02}), \qquad
h_{inv}^{21}(a)=\frac{1}{2}(3 h_{11} h_{20} +h_{02} \bar h_{02} +2 h_{11} \bar h_{11})
$$

$$
h_{inv}^{12}(a)=\frac{1}{2} (2 h_{11} ^2-h_{12} +h_{02} h_{20} +2 h_{02} \bar h_{11}+h_{11} \bar h_{20}), \qquad
h_{inv}^{03}(a)=\frac{1}{6} (3 h_{02} h_{11}-h_{03} +3 h_{02}   \bar h_{20} )
$$

$$
h_{inv}^{40}(a)=\frac{1}{24} (-15 h_{20} ^3+10 h_{20}  h_{30} -30 h_{11} h_{20} \bar h_{02} -3 h_{02}   \bar h_{02}^2
+4 h_{11} \bar h_{03}-12 h_{11} \bar h_{02} \bar h_{11})
$$

\begin{equation} \nonumber
\begin{split}
h_{inv}^{31}(a)=\frac{1}{6} (-15 h_{11}  h_{20} ^2+4 h_{11}  h_{30} -12 h_{11} ^2  \bar h_{02} +3 h_{12} \bar h_{02}-6 h_{02}  h_{20}   \bar h_{02} +h_{02} \bar h_{03} \\
-12 h_{11} h_{20} \bar h_{11} -6 h_{02} \bar h_{02} \bar h_{11} -6 h_{11} \bar h_{11}^2 +3 h_{11} \bar h_{12}-3 h_{11} \bar h_{02} \bar h_{20}) 
\end{split}
\end{equation}

\begin{equation} \nonumber
\begin{split}
h_{inv}^{22}(a)=\frac 1 4 (-12 h_{11} ^2 h_{20} +3 h_{12}  h_{20} -3 h_{02} h_{20}^2 +h_{02} h_{30} +h_{03} \bar h_{02}-9 h_{02}  h_{11} \bar h_{02}-12 h_{11}^2 \bar h_{11} \\
+4 h_{12} \bar h_{11} -6 h_{02} h_{20} \bar h_{11}-6 h_{02} \bar h_{11}^2 +2 h_{02} \bar h_{12} -3 h_{11} h_{20} \bar h_{20}-3 h_{02} \bar h_{02} \bar h_{20}-6 h_{11} \bar h_{11} \bar h_{20})
\end{split}
\end{equation}

\begin{equation} \nonumber
\begin{split}
h_{inv}^{13}(a)=\frac{1}{6} (-6 h_{11} ^3 +6 h_{11} h_{12} +h_{03} h_{20} -9 h_{02} h_{11} h_{20} -3 h_{02} ^2 \bar h_{02} +3 h_{03} \bar h_{11}-18 h_{02} h_{11} \bar h_{11} \\
-6 h_{11}^2 \bar h_{20}+3 h_{12} \bar h_{20}-3 h_{02} h_{20}   \bar h_{20} -12 h_{02} \bar h_{11} \bar h_{20}-3 h_{11}   \bar h_{20}^2 +h_{11} \bar h_{30})
\end{split}
\end{equation}

\begin{equation} \nonumber
\begin{split}
h_{inv}^{04}(a)=\frac{1}{24} (4 h_{03} h_{11} -12 h_{02} h_{11}^2 +6 h_{02} h_{12} -3 h_{02}^2 h_{20} -12 h_{02}^2 \bar h_{11} +6 h_{03} \bar h_{20}\\
-18 h_{02} h_{11} \bar h_{20} -15 h_{02} \bar h_{20}^2+4 h_{02}   \bar h_{30}) 
\end{split}
\end{equation}

\newpage

\begin{equation*}
\left|r_{2}^{40}(a)\right| = \frac{a^3 \sqrt{1 + a + \frac{4}{9 (-5 + 4 a)} + \frac{54 - 10 a (3 + a)}{9 (-9 + a (12 + (-5 + a) a))}}}{4 + a (-6 + a + a^2)}
\end{equation*}

\begin{align*}
& \left|r_{2}^{31}(a)\right| = \\
& \frac{a^3 \sqrt{4 - 6 a + 21 a^2 - 242 a^3 + 741 a^4 - 1035 a^5 + 824 a^6 - 426 a^7 + 148 a^8 - 32 a^9 + 4 a^{10}}}{(-1 + a)^2 (-5 + 4 a) \sqrt{(-1 + a) (1 + a) (-4 + a (2 + a)) (-9 + a (12 + (-5 + a) a))}}
\end{align*}

\begin{equation*}
\left|r_{2}^{22}(a)\right| =  \frac{a^3 \sqrt{1 + 2 a + 10 a^2 + 6 a^3 + 220 a^4 - 434 a^5 + 222 a^6 - 20 a^7 + a^8}}{\sqrt{-5 +4 a}(-1 + a)^3 (1 + a) (-4 + a (2 + a))}
\end{equation*}

\begin{align*}
\left|r_{2}^{13}(a)\right| & =  a^3 \big(-256 + 64 a + 1676 a^2 - 2498 a^3 - 95 a^4 + 2796 a^5 \\
& - 2219 a^6 + 187 a^7 + 730 a^8  - 550 a^9 + 200 a^{10} - 40 a^{11} + 4 a^{12}\big)^{\frac{1}{2}} \\
& \left( (-1 + a)^{\frac{5}{2}} (-5 + 4 a) (-4 + a (2 + a)) \sqrt{(1 + a) (-9 + a (12 + (-5 + a) a))}\right)^{-1}
\end{align*}

\begin{equation*}
\left|r_{2}^{04}(a)\right| = \frac{ a^3 \sqrt{-11 + a (-6 + a (27 + a (-17 + 4 a)))}}{(-4 + a (2 + a)) \sqrt{(-1 + a) (-5 + 4 a) (-9 + a (12 + (-5 + a) a))}}
\end{equation*}

\end{document}